\documentclass{article}%
\usepackage{amsmath}
\usepackage{amsfonts}
\usepackage{amssymb}
\usepackage{graphicx}%
\newtheorem{theorem}{Theorem}[section]

\newtheorem{definition}[theorem]{Definition}
\newtheorem{example}[theorem]{Example}

\newtheorem{lemma}[theorem]{Lemma}

\newtheorem{remark}[theorem]{Remark}

\newenvironment{proof}[1][Proof]{\textbf{#1.} }{\ \rule{0.5em}{0.5em}}

\begin{document}

\author{Graham Niblo\\Faculty of Mathematical Studies, University of Southampton\\Southampton SO17 1BJ, UK\\email: gan@mathematics.southampton.ac.uk
\and Michah Sageev\\Department of Mathematics, Technion\\Haifa 32000, Israel\\email: sageevm@techunix.technion.ac.il
\and Peter Scott\thanks{Partially supported by NSF grants DMS 034681 and 9626537}\\Mathematics Department, University of Michigan\\Ann Arbor, Michigan 48109, USA.\\email:pscott@umich.edu
\and Gadde A. Swarup\\Department of Mathematics and Statistics\\University of Melbourne\\Victoria 3010, Australia.\\email:gadde@ms.unimelb.edu.au }
\title{Minimal cubings}
\maketitle
\date{}

\begin{abstract}
We combine ideas of Scott and Swarup on good position for almost invariant
subsets of a group with ideas of Sageev on constructing cubings from such
sets. We construct cubings which are more canonical than in Sageev's original
construction. We also show that almost invariant sets can be chosen to be in
very good position.

\end{abstract}

Let $G$ be a finitely generated group, and let $H_{1},\ldots,H_{n}$ be
subgroups. For $i=1,\ldots,n$, let $X_{i}$ be a nontrivial $H_{i}$--almost
invariant subset of $G$. In \cite{Sageev-cubings}, Sageev gave a natural
construction of a cubing $C(X_{1},\ldots,X_{n})$ with a $G$--action which
reflects the way in which the translates of the $X_{i}$'s intersect each other.

In order to give the reader a feel for this, we start by discussing a simple
and closely related topological example. For other simple examples, the reader
is referred to Sageev's paper \cite{Sageev-cubings}. Consider a finite family
$\mathcal{F}=\{S_{1},\ldots,S_{n}\}$ of compact curves in general position on
an orientable surface $M$. There is a natural way to produce a $2$%
--dimensional cubed complex $C(\mathcal{F})$ which reflects how the $S_{i}$'s
intersect each other. Let $\widetilde{M}$ denote the universal cover of $M$,
let $\widetilde{\mathcal{F}}$ denote the pre-image of $\mathcal{F}$ in
$\widetilde{M}$, and let $D$ denote the collection of double points of the
curves in $\widetilde{\mathcal{F}}$. Then $C(\mathcal{F})$ is the dual
$2$--complex to $\widetilde{\mathcal{F}}$ in $\widetilde{M}$. This means that
$C(\mathcal{F})$ lies in $\widetilde{M}$, has one vertex in each component of
$\widetilde{M}-\widetilde{\mathcal{F}}$, and for each segment of
$\widetilde{\mathcal{F}}-D$ it has an edge which crosses this segment and no
other and joins two vertices of $C(\mathcal{F})$. Further, for each point of
the double set $D$, there is a square which contains that point and is a
$2$--cell of $C(\mathcal{F})$, and these are the only $2$--cells of
$C(\mathcal{F})$.

Now let $G$ denote $\pi_{1}(M)$. If we assume that each $S_{i}$ is essential
in $M$, then $S_{i}$ has an associated nontrivial $H$--almost invariant subset
$X_{i}$ of $G$, where $H$ equals $\pi_{1}(S_{i})$, so that $H$ is trivial or
infinite cyclic. There is a close connection between $C(\mathcal{F})$ and
Sageev's cubing $C(X_{1},\ldots,X_{n})$, although in general these cubings are
very different. Recall that $\widetilde{\mathcal{F}}$ consists of lines in
$\widetilde{M}$. Both cubings encode information about how the lines of
$\widetilde{\mathcal{F}}$ intersect. If one considers two lines of
$\widetilde{\mathcal{F}}$, the cubing $C(\mathcal{F})$ encodes very detailed
information about how they intersect, as it has a square for each double
point, but the cubing $C(X_{1},\ldots,X_{n})$ encodes only the information
about whether or not they intersect. On the other hand, if one has a family of
$k$ distinct lines in $\widetilde{\mathcal{F}}$, where $k\geq3$, and if each
line in the family meets all the others, then $C(X_{1},\ldots,X_{n})$ has a
corresponding $k$--cube, but $C(\mathcal{F})$ is always only $2$--dimensional.
However, if we assume that each component of $M-\mathcal{F}$ is not simply
connected, then the two cubings are equal. Note that this assumption implies
that no component of $\widetilde{M}-\widetilde{\mathcal{F}}$ is compact, so
that $\widetilde{\mathcal{F}}$ consists of embedded lines, and any pair of
these lines meets transversely in at most one point. Further there is no
triple of distinct lines such that each line meets the other two.

It is clear that $C(\mathcal{F})$ depends crucially on the precise
configuration of the $S_{i}$'s in $M$. For example, if the $S_{i}$'s are
disjoint, then $C(\mathcal{F})$ is $1$--dimensional, but if we homotop the
$S_{i}$'s to meet each other, then $C(\mathcal{F})$ becomes $2$--dimensional.
Thus $C(\mathcal{F})$ is not an invariant of the homotopy classes of the
curves in $\mathcal{F}$. A similar phenomenon occurs with $C(X_{1}%
,\ldots,X_{n})$. Of course, one cannot talk of almost invariant sets being
homotopic, but there is a natural idea of equivalence of almost invariant sets
which corresponds to the idea of homotopy of the $S_{i}$'s. For many groups
$G$, it is easy to give examples where $C(X_{1},\ldots,X_{n})$ is
$1$--dimensional, but if we replace each $X_{i}$ by an equivalent set $Y_{i}$,
the cubing $C(Y_{1},\ldots,Y_{n})$ is at least $2$--dimensional. Thus Sageev's
cubing depends crucially on the precise choice of the $X_{i}$'s, and is not an
invariant of the equivalence classes of the $X_{i}$'s.

In this paper, we consider the case when each of the $H_{i}$'s is finitely
generated and we show how to construct a cubing $L(X_{1},\ldots,X_{n})$ which
in most cases depends only on the equivalence classes of the $X_{i}$'s, i.e.
replacing the $X_{i}$'s by equivalent almost invariant sets yields the same
cubing. The cubing we obtain is thus more canonical than $C(X_{1},\ldots
,X_{n})$. We also show that it embeds naturally and equivariantly in
$C(X_{1},\ldots,X_{n})$ and that it is minimal in a natural sense.

Sageev's original construction depended on the partial order on the $X_{i}$'s
given by inclusion. Our construction in this paper uses Sageev's ideas but
replaces the partial order of inclusion by a partial order on the $X_{i}$'s
which is based on `almost inclusion'. Such a partial order was introduced by
Scott in \cite{Scott:TorusTheorem} in a topological context, and it played a
basic role in the purely algebraic work of Scott and Swarup in \cite{SS} and
\cite{SS-regnbhds}. In order to define this partial order, the $X_{i}$'s need
to satisfy a technical condition which Scott and Swarup called
\textquotedblleft good position\textquotedblright. In \cite{SS-regnbhds}, they
showed how to replace any finite family of almost invariant subsets of a group
by a family of equivalent almost invariant subsets which are in good position.
In this paper, we introduce an idea which we call \textquotedblleft very good
position\textquotedblright\ for almost invariant sets which is analogous to
the properties possessed by shortest curves on surfaces or by least area
surfaces in $3$--manifolds. We discuss these analogies in section
\ref{applications}. We use our new cubing to show that any finite family of
almost invariant subsets of a group can be replaced by a family of equivalent
almost invariant subsets which are in very good position. We also show how to
apply these ideas to strengthen some results of Niblo \cite{Niblo} and of
Dunwoody and Roller \cite{D-Roller}.

\section{Preliminaries}

\subsection{Almost invariant sets}

In this section, we recall the definition of an almost invariant subset of a
finitely generated group $G$, and we introduce some basic related ideas.
Throughout this paper, we will always assume that $G$ is finitely generated.
We will need several definitions which we take from \cite{SS}, but see
\cite{Scott:Intersectionnumbers} for a discussion.

\begin{definition}
Two sets $P$ and $Q$ are \textsl{almost equal} if their symmetric difference
$(P-Q)\cup(Q-P)$ is finite. We write $P\overset{a}{=}Q$.
\end{definition}

\begin{definition}
If a group $G$ acts on the right on a set $Z$, a subset $P$ of $Z$ is
\textsl{almost invariant} if $Pg\overset{a}{=}P$ for all $g$ in $G$. An almost
invariant subset $P$ of $Z$ is \textsl{nontrivial} if $P$ and its complement
$Z-P$ are both infinite. The complement $Z-P$ will be denoted simply by
$P^{\ast}$, when $Z$ is clear from the context.
\end{definition}

This idea is connected with the theory of ends of groups via the Cayley graph
$\Gamma$ of $G$ with respect to some finite generating set of $G$. (Note that
in this paper groups act on the left on covering spaces and, in particular,
$G$ acts on its Cayley graph on the left.) Using $\mathbb{Z}_{2}$ as
coefficients, we can identify $0$--cochains and $1$--cochains on $\Gamma$ with
sets of vertices or edges. A subset $P$ of $G$ represents a set of vertices of
$\Gamma$ which we also denote by $P$, and it is a beautiful fact, due to Cohen
\cite{Cohen}, that $P$ is an almost invariant subset of $G$ if and only if
$\delta P$ is finite, where $\delta$ is the coboundary operator in $\Gamma$.
Thus $G$ has a nontrivial almost invariant subset if and only if the number of
ends $e(G)$ of $G$ is at least $2$. Further $e(G)$ can be identified with the
number of nontrivial almost invariant subsets of $G$, when this count is made
correctly. If $H$ is a subgroup of $G$, we let $H\backslash G$ denote the set
of cosets $Hg$ of $H$ in $G$, i.e. the quotient of $G$ by the left action of
$H$. Of course, $G$ will no longer act on the left on this quotient, but it
will still act on the right. Thus we have the idea of an almost invariant
subset of $H\backslash G$. Further, $P$ is an almost invariant subset of
$H\backslash G$ if and only if $\delta P$ is finite, where $\delta$ is the
coboundary operator in the graph $H\backslash\Gamma$. Thus $H\backslash G$ has
a nontrivial almost invariant subset if and only if the number of ends
$e(G,H)$ of the pair $(G,H)$ is at least $2$. Considering the pre-image $X$ in
$G$ of an almost invariant subset $P$ of $H\backslash G$ leads to the
following definitions.

\begin{definition}
If $G$ is a finitely generated group and $H$ is a subgroup, then a subset $X$
of $G$ is $H$\textsl{--almost invariant} if $X$ is invariant under the left
action of $H$, and simultaneously $H\backslash X$ is an almost invariant
subset of $H\backslash G$. We may also say that $X$ is \textsl{almost
invariant over }$H$. In addition, $X$ is a \textsl{nontrivial} $H$--almost
invariant subset of $G$, if the quotient sets $H\backslash X$ and $H\backslash
X^{\ast}$ are both infinite.
\end{definition}

\begin{remark}
Note that if $X$ is a nontrivial $H$--almost invariant subset of $G$, then
$e(G,H)$ is at least $2$, as $H\backslash X$ is a nontrivial almost invariant
subset of $H\backslash G$. In fact $e(G,H)$ can be identified with the number
of nontrivial $H$--almost invariant subsets of $G$, when this count is made
correctly. See \cite{Scott-Wall:Topological} for details.
\end{remark}

\begin{definition}
If $G$ is a group and $H$ is a subgroup, then a subset $W$ of $G$ is
$H$\textsl{--finite} if it is contained in the union of finitely many left
cosets $Hg$ of $H$ in $G$.
\end{definition}

\begin{definition}
If $G$ is a group and $H$ is a subgroup, then two subsets $V$ and $W$ of $G$
are $H$\textsl{--almost equal} if their symmetric difference is $H$--finite.
\end{definition}

It will also be convenient to avoid this rather clumsy terminology sometimes,
particularly when the group $H$ is not fixed, so we make the following definition.

\begin{definition}
If $X$ is a $H$--almost invariant subset of $G\;$and $Y$ is a $K$--almost
invariant subset of $G$, and if $X$ and $Y$ are $H$--almost equal, then we
will say that $X$ and $Y$ are \textsl{equivalent} and write $X\sim Y$.
\end{definition}

\begin{remark}
Note that $H$ and $K$ must be commensurable, so that $X$ and $Y$ are also
$K$--almost equal and $(H\cap K)$--almost equal.

A more elegant and equivalent formulation is that $X$ is equivalent to $Y$ if
and only if each is contained in a bounded neighbourhood of the other. In the
context of the study of quasi-isometries, two such sets are called coarsely equivalent.

Equivalence is important because usually one is interested in an equivalence
class of almost invariant subsets of a group rather than a specific such subset.
\end{remark}

The next definitions make precise the notion of crossing of almost invariant
sets. This is an algebraic analogue of crossing of codimension--$1$ manifolds,
but it ignores \textquotedblleft inessential\textquotedblright\ crossings.

\begin{definition}
Let $X$ be an $H$--almost invariant subset of $G$ and let $Y$ be a $K$--almost
invariant subset of $G$. The four sets $X\cap Y$, $X^{\ast}\cap Y$, $X\cap
Y^{\ast}$ and $X^{\ast}\cap Y^{\ast}$ are called the \textsl{corners} of the
pair $(X,Y)$.
\end{definition}

\begin{definition}
Let $X$ be an $H$--almost invariant subset of $G$ and let $Y$ be a $K$--almost
invariant subset of $G$. We will say that $Y$ \textsl{crosses} $X$ if each of
the four corners of the pair $(X,Y)$ is not $H$--finite. Thus each of the four
corners projects to an infinite subset of $H\backslash G$.
\end{definition}

The motivation for the above definition is that when one of the four corners
is empty, we clearly have no crossing, and if one of the four corners is
``small'', then we have ``inessential crossing''. Note that $Y$ may be a
translate of $X$ in which case such crossing corresponds to the
self-intersection of a single immersion.

\begin{remark}
\label{crossingissymmetric}It is shown in \cite{Scott:Intersectionnumbers}
that if $X$ and $Y$ are nontrivial, then $X\cap Y$ is $H$--finite if and only
if it is $K$--finite. It follows that crossing of nontrivial almost invariant
subsets of $G$ is symmetric, i.e. that $X$ crosses $Y$ if and only if $Y$
crosses $X$.
\end{remark}

\begin{definition}
Let $U$ be a nontrivial $H$--almost invariant subset of $G$ and let $V$ be a
nontrivial $K$--almost invariant subset of $G$. We will say that $U\cap V$ is
\textsl{small} if it is $H$--finite.
\end{definition}

\begin{remark}
This terminology will be extremely convenient, particularly when we want to
discuss translates $U$ and $V$ of $X$ and $Y$, as we do not need to mention
the stabilisers of $U$ or of $V$. However, the terminology is symmetric in $U$
and $V$ and makes no reference to $H$ or $K$, whereas the definition is not
symmetric and does refer to $H$, so some justification is required. If $U$ is
also $H^{\prime}$--almost invariant for a subgroup $H^{\prime}$ of $G$, then
$H^{\prime}$ must be commensurable with $H$. Thus $U\cap V$ is $H$--finite if
and only if it is $H^{\prime}$--finite. In addition, Remark
\ref{crossingissymmetric} tells us that $U\cap V$ is $H$--finite if and only
if it is $K$--finite. This provides the needed justification of our terminology.

In the context of the study of quasi-isometries, the terminology
\textquotedblleft deep\textquotedblright\ is used for a subset of a metric
space which contains balls of arbitrarily large radius. One can show that
$U\cap V$ is $H$--infinite if and only if it is \textit{deep} in this sense.
\end{remark}

\subsection{Cubings\label{cubings}}

We review here the construction in \cite{Sageev-cubings}, to which the reader
is referred for details (see also \cite{NibloRoller}).

A cubed complex is a $CW$--complex formed by gluing standard Euclidean cubes
together along their faces by isometries. We further require that the boundary
of each cube is embedded in the resulting object. We do not require the
complex to be locally finite. A cubed complex is $CAT(0)$ if for every cube
$\sigma$, the link $lk(\sigma)$ of $\sigma$ satisfies the following two
conditions. There is no closed loop in $lk(\sigma)$ consisting of two edges,
and if $lk(\sigma)$ has a closed loop consisting of three edges, then this
loop bounds a triangle in $lk(\sigma)$. Finally a \textit{cubing} $C$ is a
simply connected $CAT(0)$ cubed complex. If $\sigma$ is an $n$--dimensional
cube in $C$, viewed as a standard unit cube in $\mathbb{R}^{n}$ and
$\hat{\sigma}$ denotes the barycentre of $\sigma$, then a \textit{dual cube}
in $\sigma$ is the intersection with $\sigma$ of an $(n-1)$--dimensional plane
running through $\hat{\sigma}$ and parallel to one of the $(n-1)$--dimensional
faces of $\sigma$. Given a cubing, one may consider the equivalence relation
on edges generated by the relation which declares two edges to be equivalent
if they are opposite sides of a square in $C$. Now given an equivalence class
of edges, the \textit{hyperplane} associated to this equivalence class is the
collection of dual cubes whose vertices lie on edges in the equivalence class.
It is not hard to show that hyperplanes are totally geodesic subspaces.
Moreover, in \cite{Sageev-cubings} it is shown that hyperplanes do not
self-intersect (i.e. a hyperplane meets a cube in a single dual cube) and that
a hyperplane separates a cubing into precisely two components, which we call
the \textit{half-spaces} associated to the hyperplane.

Consider a finitely generated group $G$ with subgroups $H_{1},\ldots,H_{n}$.
For $i=1,\ldots,n$, let $X_{i}$ be a nontrivial $H_{i}$--almost invariant
subset of $G$, and let $E=\{gX_{i},gX_{i}^{\ast}:g\in G,1\leq i\leq n\}$. In
\cite{Sageev-cubings}, Sageev gave a construction of a cubing from the set $E$
equipped with the partial order given by inclusion. We need the following definition.

\begin{definition}
\label{defnofultrafilter} Let $E$ be a partially ordered set, equipped with an
involution $A\rightarrow A^{\ast}$ such that $A\neq A^{\ast}$, and if $A\leq
B$ then $B^{\ast}\leq A^{\ast}$. An \textsl{ultrafilter} $V$ on $E$ is a
subset of $E$ satisfying

\begin{enumerate}
\item For every $A\in E$, we have $A\in V$ or $A^{\ast}\in V$ but not both.

\item If $A\in V$ and $A\leq B$ then $B\in V$.
\end{enumerate}
\end{definition}

Sageev constructs a cubed complex $K$ whose vertex set $K^{(0)}$ is the
collection of all ultrafilters on $E$. There is a natural action of $G$ on
$K$, and Sageev shows that a certain component $C$ of $K$ is $G$--invariant
and a cubing.

Let $K^{(0)}$ denote the collection of all ultrafilters on $E$. Construct
$K^{(1)}$ by attaching an edge to two vertices $V,V^{\prime}\in K^{(0)}$ if
and only if they differ by replacing a single element by its complement, i.e.
there exists $A\in V$ such that $V^{\prime}=(V-\{A\})\cup\{A^{\ast}\}$. Note
that the fact that $V$ and $V^{\prime}$ are both ultrafilters implies that $A$
must be a minimal element of $V$. Also if $A$ is a minimal element of $V$,
then the set $V^{\prime}=(V-\{A\})\cup\{A^{\ast}\}$ must be an ultrafilter on
$E$. Now attach $2$--dimensional cubes to $K^{(1)}$ to form $K^{(2)}$, and
inductively attach $n$--cubes to $K^{(n-1)}$ to form $K^{(n)}$. All such cubes
are attached by an isomorphism of their boundaries and, for each $n\geq2$, one
$n$--cube is attached to $K^{(n-1)}$ for each occurrence of the boundary of an
$n$--cube appearing in $K^{(n-1)}$. The complex $K$ constructed in this way
will not be connected, but one special component can be picked out in the
following way. For each element $g$ of $G$, define the ultrafilter
$V_{g}=\{A\in E:g\in A\}$. These special vertices of $K$ are called
\textit{basic}. Two basic vertices $V$ and $V^{\prime}$ of $K$ differ on only
finitely many complementary pairs of elements of $E$, so that there exist
elements $A_{1},\ldots,A_{n}$ of $E$ which lie in $V$ such that $V^{\prime}$
can be obtained from $V$ by replacing each $A_{i}$ by $A_{i}^{\ast}$. By
re-ordering the $A_{i}$'s if needed, we can arrange that $A_{1}$ is a minimal
element of $V$. It follows that $V_{1}=(V-\{A_{1}\})\cup\{A_{1}^{\ast}\}$ is
also an ultrafilter on $E$, and so is joined to $V$ by an edge of $K$. By
repeating this argument, we will find an edge path in $K$ of length $n$ which
joins $V$ and $V^{\prime}$. It follows that the basic vertices of $K$ all lie
in a single component $C$. As the collection of all basic vertices is
preserved by the action of $G$ on $K$, it follows that this action preserves
$C$. Finally, Sageev shows in \cite{Sageev-cubings} that $C$ is simply
connected and $CAT(0)$ and hence is a cubing.

At first sight, one might think that $C$ should equal $K$. To show that this
will not be the case, here are two examples.

\begin{example}
Let $E$ be the family of subsets of the integers $\mathbb{Z}$ of the form
$\{x\in\mathbb{Z}:x\leq a\}$ or $\{x\in\mathbb{Z}:x\geq b\}$, with the partial
order given by inclusion and the involution given by reflection in the
endpoint. Let $K$ and $C$ be constructed as above. Let $V$ denote the
ultrafilter on $E$ which consists of all element of $E$ of the form
$\{x\in\mathbb{Z}:x\leq a\}$. Then $V$ is not basic. In fact, $V$ differs from
any basic ultrafilter $V_{g}$ on infinitely many elements, so that $V$ is not
a vertex of $C$. Further, as $V$ has no minimal elements, it constitutes an
entire component of $K$.
\end{example}

The second example is closely related to the first, but may seem more
interesting to topologists.

\begin{example}
Let $E$ be the family of all closed half-spaces in the hyperbolic plane
$\mathbb{H}^{2}$, with the partial order given by inclusion and the involution
given by reflection of a half-space in its boundary line. Let $K$ and $C$ be
constructed as above. Let $w$ denote a point on the circle at infinity of
$\mathbb{H}^{2}$, and let $V_{w}$ denote the elements of $E$ whose closure
contains $w$. Then $V_{w}$ is not basic, and as $V_{w}$ differs from any basic
ultrafilter $V_{g}$ on infinitely many elements, it follows that $V_{w}$ is
not a vertex of $C$. Further, as $V$ has no minimal elements, it constitutes
an entire component of $K$.
\end{example}

As noted in Roller's survey article \cite{Roller}, one can characterise the
vertices of $C$ as being those ultrafilters on $E$ which satisfy the
descending chain condition. Note that the ultrafilters $V$ and $V_{w}$ in the
above two examples obviously do not satisfy the descending chain condition.

An important aspect of Sageev's construction is that one can recover the
elements of $E$ from the action of $G$ on the cubing $C$. Recall that an edge
$f$ of $C$ joins two vertices $V$ and $V^{\prime}$ if and only if there exists
$A\in V$ such that $V^{\prime}=(V-\{A\})\cup\{A^{\ast}\}$. If $f$ is oriented
towards $V^{\prime}$, we will say that $f$ \textit{exits} $A$. We let
$\mathcal{H}_{A}$ denote the hyperplane associated to the equivalence class of
$f$. This equivalence class consists of all those edges of $C$ which, when
suitably oriented, exit $A$. Now let $X$ denote an $H$--almost invariant
subset of $G$ which is an element of $E$, such that $X$ contains the identity
$e$ of $G$. Thus $X$ lies in the basic vertex $V_{e}=\{A\in E:e\in A\}$. As
$X^{\ast}$ is non-empty, it contains some element $k$ and so lies in the basic
vertex $V_{k}$. Now any path joining $V_{e}$ to $V_{k}$ must contain an edge
which exits $X$, so we can define the hyperplane $\mathcal{H}_{X}$ as above.
Let $\mathcal{H}_{X}^{+}$ denote the half-space determined by $\mathcal{H}%
_{X}$ which contains the basic vertex $V_{e}$. Recall that an edge of $C$ lies
in the equivalence class which determines $\mathcal{H}_{X}$ if and only if it
exits $X$ when suitably oriented. It follows that a vertex $V$ of $C$ lies in
$\mathcal{H}_{X}^{+}$ if and only if $X\in V$. Now we claim that the subset
$\{g\in G:gV_{e}\in\mathcal{H}_{X}^{+}\}$ of $G$ equals $X$. For
\begin{align*}
\{g  &  \in G:gV_{e}\in\mathcal{H}_{X}^{+}\}=\{g\in G:X\in gV_{e}\}=\{g\in
G:g^{-1}X\in V_{e}\}\\
&  =\{g\in G:e\in g^{-1}X\}=\{g\in G:g\in X\}=X.
\end{align*}
The following result implies that if we consider a vertex $V$ of $C$ other
than $V_{e}$, then the subset $\{g\in G:gV\in\mathcal{H}_{X}^{+}\}$ of $G$ is
still $H$--almost invariant, and although it need not be equal to $X$, it is
still equivalent to $X$.

\begin{lemma}
\label{lem:ConstructinAISets} Suppose that $G$ is a finitely generated group
which acts on a cubing $C$. Let $\mathcal{H}$ be a hyperplane in $C$ with
stabilizer $H$, let $\mathcal{H}^{+}$ and $\mathcal{H}^{-}$ denote the two
half-spaces defined by $\mathcal{H}$, and suppose that $H$ preserves each of
$\mathcal{H}^{+}$ and $\mathcal{H}^{-}$. Then, for any vertex $v$, the set
$X_{v}=\{g\in G|gv\in\mathcal{H}^{+}\}$ is almost invariant over $H$ and all
these subsets of $G$ are equivalent.
\end{lemma}

\begin{proof}
We need to show that $hX_{v}=X_{v}$, for all $h$ in $H$, and that $X_{v}a$ and
$X_{v}$ are $H$--almost equal for all $a$ in $G$.

As $H$ stabilises $\mathcal{H}^{+}$, it follows immediately that $hX_{v}%
=X_{v}$, for all $h$ in $H$.

Next consider $X_{v}-X_{v}a$. From the definition of $X_{v}$, we have that%

\[
X_{v}a=\{ga\in G|gv\in\mathcal{H}^{+}\}=\{g^{\prime}\in G|g^{\prime}a^{-1}%
v\in\mathcal{H}^{+}\}.
\]

Hence%
\[
X_{v}-X_{v}a=\{g\in G|gv\in\mathcal{H}^{+}\ \mathrm{and}\ ga^{-1}%
v\notin\mathcal{H}^{+}\}.
\]

Thus $g\in X_{v}-X_{v}a$ if and only if $\mathcal{H}$ separates $gv$ from
$ga^{-1}v$. Now there are only finitely many hyperplanes in $C$ which separate
$v$ from $a^{-1}v$. We denote these hyperplanes by $\mathcal{H}_{1}%
,\ldots,\mathcal{H}_{n}$. It follows that if $g\in X_{v}-X_{v}a$, then
$\mathcal{H}=g\mathcal{H}_{i}$, for some $i$. For any two elements $g$ and
$g^{\prime}$ such that $\mathcal{H}=g\mathcal{H}_{i}$ and $\mathcal{H}%
=g^{\prime}\mathcal{H}_{i}$, we have that $g^{\prime}g^{-1}\mathcal{H}%
=\mathcal{H}$, so that $g^{\prime}g^{-1}\in H$ and $Hg=Hg^{\prime}$. It
follows that $X_{v}-X_{v}a$ is contained in $HF$ for some finite set $F$, and
so is $H$--finite. Similarly, $X_{v}a-X_{v}$ is $H$--finite. As this holds for
any element $a$ of $G$, it follows that $X_{v}$ is almost invariant over $H$,
as required.

Now let $v$ and $w$ denote two vertices of $C$, and let $k$ be an element of
$X_{v}-X_{w}$. Thus $kv\in\mathcal{H}^{+}$ and $kw\notin\mathcal{H}^{+}$.
Hence $\mathcal{H}$ separates $kv$ from $kw$, so that $k^{-1}\mathcal{H}$
separates $v$ and $w$. As in the above argument, it follows that $X_{v}-X_{w}$
is $H$--finite. Similarly, $X_{w}-X_{v}$ is $H$--finite. It follows that
$X_{v}$ and $X_{w}$ are equivalent, which completes the proof of the lemma.
\end{proof}

\section{The new partial order\label{newpartialorder}}

In this section, we recall some of the ideas of Scott and Swarup in
\cite{SS-regnbhds} and \cite{SS}.

Consider a finitely generated group $G$ with finitely generated subgroups
$H_{1},\ldots,H_{n}$. For $i=1,\ldots,n$, let $X_{i}$ be a nontrivial $H_{i}%
$--almost invariant subset of $G$, and let $E=\{gX_{i},gX_{i}^{\ast}:g\in
G,1\leq i\leq n\}$. As $E$ is a collection of subsets of $G$, it has a natural
partial order induced by inclusion. But one can sometimes define a more
interesting partial order. The idea is to define $U\leq V$ when $U$ is
\textquotedblleft nearly\textquotedblright\ contained in $V$. Precisely, we
want $U\leq V$ if $U\cap V^{\ast}$ is small. However, an obvious difficulty
arises when two of the corners $U^{(\ast)}\cap V^{(\ast)}$ are small, as we
have no way of deciding between two possible inequalities. It turns out that
we can avoid this difficulty if we know that whenever two of the corners of
$U\ $and $V$ are small, then one of them is empty. Thus we consider the
following condition on $E$:

Condition (*): If $U$ and $V$ are in $E$, and two of their corners are small,
then one of their corners is empty.

If $E$ satisfies Condition (*), we will say that the family $X_{1}\ldots
,X_{n}$ is in \textit{good position}.

Assuming that this condition holds, we can define a relation $\leq$ on $E$ by
saying that $U\leq V$ if and only if $U\cap V^{\ast}$ is empty or is the only
small set among the four corners of $U$ and $V$. Despite the seemingly
artificial nature of this definition, one can show that $\leq$ is a partial
order on $E$. This is not entirely trivial, but the proof is in Lemma 1.14 of
\cite{SS}. Condition (*) plays a key role in the proof. If $U\leq V$ and
$V\leq U$, it is easy to see that we must have $U=V$, using the fact that
$E\;$satisfies Condition (*). Most of the proof of Lemma 1.14 of \cite{SS} is
devoted to showing that $\leq$ is transitive.

We will need the following fact, which follows immediately from Lemma 2.31 of
\cite{SS-regnbhds}. Note that the number $D$ is independent of the element $g$
of $G$.

\begin{lemma}
\label{containedinboundednbhdofA}Let $G$ be a finitely generated group with
finitely generated subgroups $H$ and $K$, a nontrivial $H$--almost invariant
subset $A$ and a nontrivial $K$--almost invariant subset $U$. Let $\Gamma$
denote the Cayley graph of $G$ with respect to some finite generating set.
Then there is $D>0$, such that if $gU\leq A$, then $gU$ is contained in the
$D$--neighbourhood of $A$ in $\Gamma$.
\end{lemma}

\begin{remark}
This result will play a key role in our construction of a cubing in section
\ref{CubingConstruction}. This explains why we need to restrict our attention
to almost invariant subsets of $G$ which are over finitely generated subgroups.
\end{remark}

In general, the family $X_{1},\ldots,X_{n}$ need not be in good position, but
we will use the results in \cite{SS} to show that we can find almost invariant
sets $Y_{1},\ldots,Y_{n}$ such that $Y_{i}$ is equivalent to $X_{i}$ and the
$Y_{i}$'s are in good position. We will also show that the partial order
obtained is unique in most cases. Scott and Swarup did not state such results
in \cite{SS}, as they were concentrating on almost invariant sets associated
to splittings, but all the arguments needed are essentially there.

It turns out that the case when $n=1$ contains almost all of the difficulties,
so we will start by discussing that case. Let $H$ be a finitely generated
subgroup of $G$, and let $X$ be a $H$--almost invariant subset of $G$. If $X$
is not in good position, there must be two translates $U\;$and $V$ of $X$ such
that two of their corners are small, and neither is empty. If $U\cap V$ is one
of the two small corners, the other must be $U^{\ast}\cap V^{\ast}$, as
otherwise $U$ or $V$ would be small which contradicts the fact that $X$ is
nontrivial. Similarly, if $U\cap V^{\ast}$ is one of the two small corners,
the other must be $U^{\ast}\cap V$. It follows that $U$ is equivalent to $V$
or to $V^{\ast}$. This naturally leads one to consider the subgroup
$\mathcal{K}(X)$ of $G$ defined by $\mathcal{K}(X)=\{g\in G:gX\sim X$
\textrm{or} $X^{\ast}\}$. It will also be convenient to consider the subgroup
$\mathcal{K}_{0}(X)=\{g\in G:gX\sim X\}$ of $\mathcal{K}$, so that the index
of $\mathcal{K}_{0}$ in $\mathcal{K}$ is at most $2$. We will say that the
collection $E(X)$ of all translates of $X$ and $X^{\ast}$ is \textit{nested
with respect to} $\mathcal{K}$, if for any $k\in\mathcal{K}$, one of the four
corners of $X$ and $kX$ is empty. It is clear that $X$ is in good position if
and only if $E(X)$ is nested with respect to $\mathcal{K}$. The following
lemma summarises results proved by Scott and Swarup in the proof of
Proposition 2.14 of \cite{SS}.

\begin{lemma}
\label{Prop2.14ofSS}(Scott-Swarup) Let $G$ be a finitely generated group with
a finitely generated subgroup $H$, and let $X$ be a nontrivial $H$--almost
invariant subset of $G$.

\begin{enumerate}
\item If $H$ has finite index in $\mathcal{K}$, there is an almost invariant
subset $W$ of $G$ with stabiliser $\mathcal{K}_{0}$ which is equivalent to
$X$, such that $E(W)$ is nested with respect to $\mathcal{K}$.

\item If $H$ has infinite index in $\mathcal{K}$, then $\mathcal{K}$ has
finite index in $G$, and there is a subgroup $H^{\prime\prime}$ of
$\mathcal{K}$ which is commensurable with $H$ and normal in $\mathcal{K}.$
Further, $H^{\prime\prime}\backslash\mathcal{K}$ is isomorphic to $\mathbb{Z}$
or to $\mathbb{Z}_{2}\ast\mathbb{Z}_{2}$. In the first case $\mathcal{K}%
=\mathcal{K}_{0}$, and in the second case $\mathcal{K}_{0}$ has index $2$ in
$\mathcal{K}$. There is an almost invariant subset $W$ of $G$ with stabiliser
$H^{\prime\prime}$ which is equivalent to $X$, such that $E(W)$ is nested with
respect to $\mathcal{K}$.
\end{enumerate}
\end{lemma}

Now we can prove the following result.

\begin{lemma}
\label{canobtaingoodpositionifn=1}Let $G$ be a finitely generated group with a
finitely generated subgroup $H$, and let $X$ be a nontrivial $H$--almost
invariant subset of $G$. Then $X$ is equivalent to a $K$--almost invariant
subset $W$ of $G$ which is in good position. Thus the set $E(W)$ of all
translates of $W$ and $W^{\ast}$ has the partial order $\leq$ described above.
\end{lemma}

\begin{proof}
Lemma \ref{Prop2.14ofSS} shows that in all cases, there is an almost invariant
subset $W$ of $G$ which is equivalent to $X$ such that $E(W)$ is nested with
respect to $\mathcal{K}(X)$. As $X$ and $W$ are equivalent, the subgroups
$\mathcal{K}(X)$ and $\mathcal{K}(W)$ are equal, so that $E(W)$ is nested with
respect to $\mathcal{K}(W)$. As remarked just before the statement of Lemma
\ref{Prop2.14ofSS}, this implies that $W$ is in good position, which completes
the proof.
\end{proof}

We would like to show that the partial order obtained by applying the above
result is unique. More precisely, if $Y\ $and $Z$ are equivalent to $X$ and in
good position, we want to show that there is a $G$--equivariant bijection
between $E(Y)\ $and $E(Z)$ which preserves complementation and the partial
orders. It is natural to attempt to define such a map $\varphi:E(Y)\rightarrow
E(Z)$, by sending $Y$ to $Z$, and extending appropriately. If it is to be
$G$--equivariant, it must send $gY$ to $gZ$ for every $g$ in $G$. This
immediately raises a potential problem, which is that it seems possible that
$gY=Y$, but $gZ\neq Z$. However the following result shows that this cannot occur.

\begin{lemma}
\label{equivalentgoodpositionimpliessamestabiliser}Let $G$ be a finitely
generated group, let $Y$ and $Z$ be equivalent almost invariant subsets of $G$
each of which is in good position. Then the stabilisers of $Y\ $and $Z$ are equal.
\end{lemma}

\begin{proof}
Let $K$ and $L$ denote the stabilisers of $Y\ $and $Z$ respectively, so that
$K$ and $L$ must be commensurable subgroups of $G$. Let $k$ denote an element
of $K$, so that $kY=Y$. As $Z$ is equivalent to $Y$, it follows that $kZ$ is
equivalent to $Z$. As $Z$ is in good position, we must have $kZ=Z$, or
$Z\subset kZ$ or $kZ\subset Z$. As $K\ $and $L$ are commensurable, some power
$k^{n}$ of $K$ must lie in $L$, so that $k^{n}Z=Z$. It follows that in all
cases we must have $kZ=Z$, so that $k$ lies in $L$. Thus $K$ is contained in
$L$. Similarly $L$ is contained in $K$, so that $K=L$ as required.
\end{proof}

Now we return to the question of the uniqueness of the partial order on $E(W)$
obtained by applying Lemma \ref{canobtaingoodpositionifn=1}. Suppose that
$Y\ $and $Z$ are equivalent to $X$ and in good position. We want to define a
bijection $\varphi:E(Y)\rightarrow E(Z)$, which is $G$--equivariant and
preserves complementation. If $\varphi$ sends $Y$ to $Z$ it must also send
$gY$ to $gZ$ and $gY^{\ast}$ to $gZ^{\ast}$, for every $g$ in $G$. The fact
that the stabilisers of $Y\ $and $Z$ are equal implies that this gives a well
defined map on the translates of $Y$. There is still a potential problem,
which is that it seems possible that $gY=Y^{\ast}$, but $gZ\neq Z^{\ast}$. If
this does not occur, it is clear that we do have a well defined map from
$E(Y)$ to $E(Z)$ which is $G$--equivariant and preserves complementation. In
order to discuss the general situation, we will use the following piece of
terminology which Scott and Swarup introduced in \cite{SS-regnbhds}.

\begin{definition}
\label{defnofinvertible}If $X$ is an $H$--almost invariant subset of a group
$G$, then $X$ is \textsl{invertible} if there is an element $g$ in $G$ such
that $gX=X^{\ast}$.
\end{definition}

Note that in \cite{SS-regnbhds}, Scott and Swarup only used this term when $X$
was associated to a splitting, but in this paper, we will not make that restriction.

Our previous discussion shows that if $Y$ is not invertible, then we have a
well defined map $\varphi:E(Y)\rightarrow E(Z)$, described by sending $gY$ to
$gZ$ and $gY^{\ast}$ to $gZ^{\ast}$, for every $g$ in $G$. If, in addition,
$Z$ is not invertible, then the same comment applies to the inverse map
showing that $\varphi$ must be a bijection, which is $G$--equivariant and
preserves complementation. It is also clear that $\varphi(U)$ is equivalent to
$U$ for every $U$ in $E(Y)$. We will say that $\varphi$ \textit{preserves
equivalence classes}.

Now we can prove our first uniqueness result for partial orders.

\begin{lemma}
\label{goodpositionyieldsuniquepoifn=1}Let $G$ be a finitely generated group
with a finitely generated subgroup $H$. Let $X$ be a nontrivial $H$--almost
invariant subset of $G$, and suppose that $X$ is equivalent to $Y$ and to $Z$
such that each of $Y$ and $Z$ is in good position. In addition, suppose that
$Y$ and $Z$ are both not invertible. Then there is a $G$--equivariant
bijection $\varphi:E(Y)\rightarrow E(Z)$ which preserves the partial order
$\leq$ and preserves complementation and equivalence classes.
\end{lemma}

\begin{proof}
As discussed above, we can define a $G$--equivariant bijection $\varphi
:E(Y)\rightarrow E(Z)$, by sending $gY$ to $gZ$ and $gY^{\ast}$ to $gZ^{\ast}$
for every $g$ in $G$, and $\varphi$ also preserves complementation and
equivalence classes. In many cases, $\varphi$ is already order preserving, but
if it is not we will describe a simple modification of $\varphi$ which will
arrange this.

Let $U$ and $V$ denote elements of $E(Y)$. As $U$ is equivalent to
$\varphi(U)$ and $V$ is equivalent to $\varphi(V)$, it follows that a corner
of $U\ $and $V$ is small if and only if the corresponding corner of
$\varphi(U)$ and $\varphi(V)$ is small. Hence $U$ and $V$ are comparable in
$E(Y)$ if and only if $\varphi(U)$ and $\varphi(V)$ are comparable in $E(Z)$.
Further, it follows that $\varphi$ is order preserving, except possibly when
there are $U\;$and $V$ such that two of the four corners of $U$ and $V$ are
small. If this happens, then $U\;$and $V$ must be equivalent, and we again
consider the group $\mathcal{K}(X)=\{g\in G:gX\sim X$ \textrm{or} $X^{\ast}%
\}$. Note that as $X$, $Y$ and $Z\ $are equivalent, the groups $\mathcal{K}%
(X)$, $\mathcal{K}(Y)$ and $\mathcal{K}(Z)$ are all equal. We denote this
group by $\mathcal{K}$. We also have the subgroup $\mathcal{K}_{0}=\{g\in
G:gX\sim X\}$ of $\mathcal{K}$, whose index in $\mathcal{K}$ is at most $2$.

Suppose that $H$ has finite index in $\mathcal{K}$. Then part 1) of Lemma
\ref{Prop2.14ofSS} implies that there is an almost invariant subset $W$ of $G$
with stabiliser $\mathcal{K}_{0}$ which is equivalent to $X$ and in good
position. The fact that $W$ is in good position combined with Lemma
\ref{equivalentgoodpositionimpliessamestabiliser} implies that the stabilisers
of $Y$ and $Z$ also equal $\mathcal{K}_{0}$. If $\mathcal{K}=\mathcal{K}_{0}$,
it follows that $\varphi$ is order preserving, because there are no distinct
equivalent elements of $E(Y)$. If $\mathcal{K}_{0}$ has index $2$ in
$\mathcal{K}$, it is possible that $\varphi$ is not order preserving, so we
need some special arguments. If $k$ denotes an element of $\mathcal{K}%
-\mathcal{K}_{0}$, then $kY^{\ast}$ is equivalent to $Y$. As $Y$ is in good
position, we must have $kY^{\ast}\subset Y$ or $Y\subset kY^{\ast}$. Note that
as we are assuming that $Y$ is not invertible, we cannot have $Y=kY^{\ast}$.
We can suppose that $kY^{\ast}\subset Y$, by replacing $k$ by $k^{-1}$ and $Y$
by $Y^{\ast}$, if necessary. Thus either $\varphi$ is order preserving, or
this fails to hold only in that $kY^{\ast}\subset Y$ but $Z\subset kZ^{\ast}$,
for all $k\in\mathcal{K}-\mathcal{K}_{0}$. If $\varphi$ is not order
preserving, we replace $Z$ by $Z^{\prime}=kZ$ and we replace $Y$ by
$Y^{\prime}=Y^{\ast}.$ As $Y^{\prime}$ and $Z^{\prime}$ are each in good
position, and equivalent to each other, there is a natural $G$--equivariant
bijection $\varphi^{\prime}:E(Y^{\prime})\rightarrow E(Z^{\prime})$ sending
$Y^{\prime}$ to $Z^{\prime}$ which must be order preserving, except possibly
when one compares $Y^{\prime}$, $kY^{\prime}$ and $Z^{\prime}$, $kZ^{\prime}$,
where $k\in\mathcal{K}-\mathcal{K}_{0}.$ Now the inclusion $kY^{\ast}\subset
Y$ tells us that $kY^{\prime}\subset(Y^{\prime})^{\ast}$, and the inclusion
$Z\subset kZ^{\ast}$ tells us that $kZ^{\prime}=k^{2}Z=Z\subset kZ^{\ast
}=\left(  Z^{\prime}\right)  ^{\ast}$. We conclude that $\varphi^{\prime}$ is
order preserving, and preserves complementation and equivalence classes.

Now suppose that $H$ has infinite index in $\mathcal{K}$. Then part 2) of
Lemma \ref{Prop2.14ofSS} tells us that $\mathcal{K}$ has finite index in $G$,
and there is a subgroup $H^{\prime\prime}$ of $\mathcal{K}$ which is
commensurable with $H$ and normal in $\mathcal{K}.$ Further, $H^{\prime\prime
}\backslash\mathcal{K}$ is isomorphic to $\mathbb{Z}$ or to $\mathbb{Z}%
_{2}\ast\mathbb{Z}_{2}$. In the first case $\mathcal{K}=\mathcal{K}_{0}$, and
in the second case $\mathcal{K}_{0}$ has index $2$ in $\mathcal{K}$. It also
implies that there is an almost invariant subset $W$ of $G$ with stabiliser
$H^{\prime\prime}$ which is equivalent to $X$ and in good position. As before,
it follows that the stabilisers of $Y$ and $Z$ must also equal $H^{\prime
\prime}$. The facts that $H^{\prime\prime}$ is normal in $\mathcal{K}$ with
quotient a group with two ends, and that $\mathcal{K}$ has finite index in
$G$, imply that $e(G,H^{\prime\prime})=2$. If $\mathcal{K}=\mathcal{K}_{0}$,
we let $\lambda$ denote an element of $\mathcal{K}$ which maps to a generator
of $H^{\prime\prime}\backslash\mathcal{K}$, and we choose $\lambda$ so that
$Y\subset\lambda Y$. Then either $Z\subset\lambda Z$ or $\lambda Z\subset Z$.
As $Y$ and $Z$ are equivalent, there is a number $D$ such that $Y$ and $Z$
each lie in the $D$--neighbourhood of the other. Hence the unions
$\bigcup_{n\geq1}${$\lambda$}$^{n}Y$ and $\bigcup_{n\geq1}\lambda^{n}Z$ each
lie in the $D$--neighbourhood of the other. As $Y\subset\lambda Y$, and
$e(G,H^{\prime\prime})=2$, the union $\bigcup_{n\geq1}\lambda^{n}Y$ equals
$G$. It follows that the union $\bigcup_{n\geq1}\lambda^{n}Z$ also equals $G$,
so that the inclusion $\lambda Z\subset Z$ is impossible. Thus $Z\subset
\lambda Z$, which implies that $\varphi$ is order preserving.

If $\mathcal{K}\neq\mathcal{K}_{0}$, so that $H^{\prime\prime}\backslash
\mathcal{K}$ is $\mathbb{Z}_{2}\ast\mathbb{Z}_{2}$, the situation is more
complicated. Fix an element $k$ of $\mathcal{K}-\mathcal{K}_{0}$, so that $kY$
is equivalent to $Y^{\ast}$. As $Y$ is in good position, we must have either
$Y\subset kY^{\ast}$ or $Y^{\ast}\subset kY$. (Again the assumption that $Y$
is not invertible implies that we cannot have $Y=kY^{\ast}$.) Similarly, for
each integer $n$, we must have $\lambda^{n}Y\subset k\lambda^{n}Y^{\ast}$ or
$\lambda^{n}Y^{\ast}\subset k\lambda^{n}Y$. Suppose that $\lambda^{n}Y\subset
k\lambda^{n}Y^{\ast}$, for some $n$. As $Y\subset\lambda Y$, it follows that
$\lambda^{n}Y\subset\lambda^{m+n}Y$, for every $m\geq1$, and so $k\lambda
^{n}Y\subset k\lambda^{m+n}Y$, for every $m\geq1$. As the union of the
$\lambda^{m+n}Y$, for $m\geq1$, equals $G$, so does the union of the
$k\lambda^{m+n}Y$, for $m\geq1$. It follows that we cannot have $\lambda
^{m+n}Y\subset k\lambda^{m+n}Y^{\ast}$, for every $m\geq1$. In particular, the
inclusion $\lambda^{n}Y\subset k\lambda^{n}Y^{\ast}$ cannot hold for all
values of $n$. Similarly, the inclusion $\lambda^{n}Y^{\ast}\subset
k\lambda^{n}Y$ cannot hold for all values of $n$. If $\lambda^{N}Y\subset
k\lambda^{N}Y^{\ast}$ for some integer $N$, then $\lambda^{n}Y\subset
k\lambda^{n}Y^{\ast}$ whenever $n\leq N$. It follows that there is an integer
$N(Y)$ such that $\lambda^{n}Y\subset k\lambda^{n}Y^{\ast}$ whenever $n\leq
N(Y)$, and $\lambda^{n}Y^{\ast}\subset k\lambda^{n}Y$ whenever $n>N(Y)$. A
similar discussion for $Z$ yields an integer $N(Z)$ such that $\lambda
^{n}Z\subset k\lambda^{n}Z^{\ast}$ whenever $n\leq N(Z)$, and $\lambda
^{n}Z^{\ast}\subset k\lambda^{n}Z$ whenever $n>N(Z)$. If $N(Y)=N(Z)$, it is
now easy to see that $\varphi$ is order preserving. Otherwise, we let $d$
denote $N(Z)-N(Y)$ and let $Z^{\prime}$ denote $\lambda^{d}Z$, so that
$Z^{\prime}$ is equivalent to $Z$, and let $\varphi^{\prime}:E(Y)\rightarrow
E(Z^{\prime})$ be the equivariant bijection which sends $Y$ to $Z^{\prime}$.
As $N(Z^{\prime})=N(Y)$, it follows that $\varphi^{\prime}$ is order
preserving, and so is the required order preserving bijection from $E(Y)$ to
$E(Z)$.
\end{proof}

The above result shows that when one replaces $X$ by an almost invariant set
in good position, one obtains a unique partial order if we do not allow
invertible almost invariant sets. We now discuss the general situation.
Clearly if $Y\ $and $Z$ are equivalent to $X$ and one is invertible and the
other is not, we do not obtain exactly the same partial order, so we now
restrict attention to the case where both $Y\ $and $Z$ are invertible.

\begin{lemma}
\label{goodpositionyieldsalmostuniquepoifn=1invertible}Let $G$ be a finitely
generated group with a finitely generated subgroup $H$. Let $X$ be a
nontrivial $H$--almost invariant subset of $G$, and suppose that $X$ is
equivalent to $Y$, $Z$ and $V$ such that each of $Y$, $Z$ and $V$ is in good
position. Thus $\leq$ defines a partial order on $E(Y)$, $E(Z)$ and $E(V)$. In
addition, suppose that $Y$, $Z$ and $V$ are each invertible. Then one of the
following holds:

\begin{enumerate}
\item There are $G$--equivariant bijections between $E(Y)$, $E(Z)$ and $E(V)$
which preserve complementation, ordering and equivalence classes.

\item $H$ has infinite index in $\mathcal{K}$ and there is a $G$--equivariant
bijection between two of $E(Y)$, $E(Z)$ and $E(V)$ which preserves
complementation, ordering and equivalence classes.
\end{enumerate}
\end{lemma}

\begin{remark}
This means that in case 1) there is only one partially ordered set as in Lemma
\ref{goodpositionyieldsuniquepoifn=1}, and in case 2) there are at most two
possible partially ordered sets. The case of two distinct partial orders can
occur. The simplest example occurs when $G$ is $\mathbb{Z}_{2}\ast
\mathbb{Z}_{2}$ and $H$ is trivial.
\end{remark}

\begin{proof}
For simplicity, we start by considering $Y\ $and $Z$ only. The assumption that
$Y$ and $Z$ are both invertible implies that $\mathcal{K}_{0}$ has index $2$
in $\mathcal{K}$. It is no longer obvious that we can define a $G$%
--equivariant map $\varphi:E(Y)\rightarrow E(Z)$, by sending $gY$ to $gZ$ and
$gY^{\ast}$ to $gZ^{\ast}$ for every $g$ in $G$, because it is possible that
there is $g$ in $G$ such that $gY=Y^{\ast}$ but $gZ\neq Z^{\ast}$.

If $H$ has finite index in $\mathcal{K}$, then as in the proof of Lemma
\ref{goodpositionyieldsuniquepo} the stabilisers of $Y$ and $Z$ must both
equal $\mathcal{K}_{0}$. As each of $Y$ and $Z$ is invertible, it follows that
$kY=Y^{\ast}$ and $kZ=Z^{\ast}$ for every $k$ in $K-K_{0}$. Hence $\varphi$
can be defined as above, and it is a $G$--equivariant bijection. It is also
order preserving because there are no distinct equivalent elements of $E(Y)$.

Now suppose that $H$ has infinite index in $\mathcal{K}$. Then as in the proof
of Lemma \ref{goodpositionyieldsuniquepo} the stabilisers of $Y$ and $Z$ must
equal $H^{\prime\prime}$. In this case, it is possible that $\varphi$ cannot
be defined as above, because the elements which invert $Y$ and $Z$ need not be
the same. As in the case when $Y\ $and $Z$ were not invertible, we let
$\lambda$ denote an element of $\mathcal{K}$ which maps to a generator of
$H^{\prime\prime}\backslash\mathcal{K}$, and we choose $\lambda$ so that
$Y\subset\lambda Y$. As in that case, it follows that $Z\subset\lambda Z$. Now
let $k$ denote an element of $\mathcal{K}-\mathcal{K}_{0}$ such that
$kY=Y^{\ast}$. As $Y\subset\lambda Y$ and so $\lambda Y^{\ast}\subset Y^{\ast
}$, it is clear that $k$ cannot invert $\lambda^{n}Y$, for any $n\neq0$. If
$kZ=Z^{\ast}$, then $\varphi$ can be defined as above and is a $G$%
--equivariant bijection. Further it is easy to see that $\varphi$ is order
preserving. If $kZ\neq Z^{\ast}$, the fact that $Z$ is invertible means that
there is an integer $n\neq0$ such that $k\lambda^{n}Z=Z^{\ast}$. If $n$ is
even, say $n=2m$, this is equivalent to the equation $k\lambda^{m}%
Z=\lambda^{m}Z^{\ast}$, and we let $Z^{\prime}=\lambda^{m}Z$. We can now
define $\varphi^{\prime}:E(Y)\rightarrow E(Z^{\prime})$ to send $gY$ to
$gZ^{\prime}$ and $gY^{\ast}$ to $gZ^{\prime\ast}$, and $\varphi^{\prime}$ is
a $G$--equivariant bijection which preserves complementation and is order
preserving. As $E(Z^{\prime})=E(Z)$, this is the required bijection. However,
if $n$ is odd, this cannot be done.

To complete the proof of the lemma, we consider all three of $Y$, $Z$ and $V$.
If $H$ has finite index in $\mathcal{K}$, the above proof applies to each pair
to show that the required $G$--equivariant bijections exist. If $H$ has
infinite index in $\mathcal{K}$, we consider the preceding paragraph. Choose
$\lambda$ and $k$ as described there. There is an integer $n$ such that
$k\lambda^{n}Z=Z^{\ast}$. Similarly, there is an integer $r$ such that
$k\lambda^{r}V=V^{\ast}$. If either of $n$ or $r$ is even, the preceding
paragraph provides a $G$--equivariant bijection between $E(Y)$ and one of
$E(Z)$ or $E(V)$. If both $n$ and $r$ are odd, we let $k^{\prime}$ denote
$k\lambda$, so that we have the equations $k^{\prime}\lambda^{n-1}Z=Z^{\ast}$
and $k^{\prime}\lambda^{r-1}Z=Z^{\ast}$. As $n-1$ and $r-1$ are both even, say
$n-1=2m$ and $r-1=2s$, we let $Z^{\prime}=\lambda^{m}Z$ and $V^{\prime
}=\lambda^{s}V$. Thus $k^{\prime}$ inverts $Z^{\prime}$ and inverts
$V^{\prime}$. Now we can define $\varphi^{\prime}:E(Z^{\prime})\rightarrow
E(V^{\prime})$ to send $gZ^{\prime}$ to $gV^{\prime}$ and $gZ^{\prime\ast}$ to
$gV^{\prime\ast}$, and $\varphi^{\prime}$ is the required $G$--equivariant
bijection $E(Z)=E(V)$.
\end{proof}

The above discussion shows that if one considers all possible ways of
replacing $X$ by an almost invariant set in good position, only one partially
ordered set can be obtained in this way, unless $X$ is equivalent to an
invertible almost invariant set. In this case, at most two partially ordered
sets can be obtained with $Y$ invertible. Thus in all cases, at most three
partially ordered sets can be obtained by replacing $X$ by an almost invariant
set in good position.

This completes our discussion of good position when one starts with a single
almost invariant subset of $G$. It is now easy to extend this to the general case.

\begin{lemma}
\label{canobtaingoodposition}Let $G$ be a finitely generated group with
finitely generated subgroups $H_{1},\ldots,H_{n}$. For $i=1,\ldots,n$, let
$X_{i}$ be a nontrivial $H_{i}$--almost invariant subset of $G$. Then each
$X_{i}$ is equivalent to a $K_{i}$--almost invariant subset $Y_{i}$ of $G$
such that the $Y_{i}$'s are in good position. Thus the set $E(Y_{1}%
,\ldots,Y_{n})$ of all translates of all the $Y_{i}$'s and their complements
has the partial order $\leq$ described above.
\end{lemma}

\begin{proof}
By Lemma \ref{canobtaingoodpositionifn=1}, we can replace each $X_{i}$ by an
equivalent almost invariant set $Y_{i}$, such that each $Y_{i}$ is in good
position. Thus for each $i$, the set $E(Y_{i})$ of all translates of $Y_{i}$
and $Y_{i}^{\ast}$ satisfies Condition (*). Suppose that the set
$E(Y_{1},\ldots,Y_{n})$ of all translates of all the $Y_{i}$'s and
$Y_{i}^{\ast}$'s does not satisfy Condition (*). Then there exist distinct $i$
and $j$ and translates $U\ $and $V$ of $Y_{i}$ and $Y_{j}$ respectively such
that two of their corners are small, and neither is empty. As before, this
implies that $U$ is equivalent to $V$ or to $V^{\ast}$, so that $Y_{i}$ is
equivalent to some translate of $Y_{j}$ or $Y_{j}^{\ast}$. In this case we
simply replace $Y_{i}$ by the same translate of $Y_{j}$ or $Y_{j}^{\ast}$. By
repeating this process, we will be able to arrange that the collection
$Y_{1},\ldots,Y_{n}$ is also in good position, as required.
\end{proof}

In the preceding proof, it may seem that we took the easy way out by simply
replacing $Y_{i}$ by a translate of $Y_{j}$ or $Y_{j}^{\ast}$. However the
following simple example shows that there are cases when there is no other way
to arrange that the $Y_{i}$'s are in good position.

\begin{example}
Let $G$ denote the integers under addition and let $H$ denote the trivial
subgroup of $G$. As $G$ has two ends, it has nontrivial almost invariant
subsets over $H$. The natural examples are sets of the form $L_{a}=\{n\in
G:n\leq a\}$ or $R_{a}=\{n\in G:n\geq a\}$ for some integer $a$. If $X$ is an
almost invariant subset of $G$ over $H$ which is in good position, it is easy
to see that $X$ must be one of the sets $L_{a}$ or $R_{a}$, for some $a$. Thus
the set $E(X)$ of all translates of $X$ and $X^{\ast}$ consists of all the
sets $L_{a}$ and $R_{a}$. It follows that it is impossible to have two almost
invariant subsets $X_{1}$ and $X_{2}$ of $G$ such that $E(X_{1},X_{2})$
satisfies Condition (*) unless $X_{2}$ is some translate of $X_{1}$ or
$X_{1}^{\ast}$. Thus in this group there is simply not room for more than one
almost invariant set to be in good position.
\end{example}

The above example suggests that if we want the $Y_{i}$'s we choose in Lemma
\ref{canobtaingoodposition} to be in good position and to reflect the
properties of the $X_{i}$'s, then we should exclude the possibility that there
are $X_{i}$ and $X_{j}$, with $i\neq j$, such that some translate of $X_{i}$
is equivalent to $X_{j}$ or $X_{j}^{\ast}$. If this occurs, we will say that
the $G$--orbits of $X_{i}$ and $X_{j}$ are \textit{parallel}. We use this word
because we are thinking of parallel $G$--orbits as corresponding to homotopic
curves on a surface. The following simple uniqueness result covers most
situations. However, if one allows some of the $Y_{i}$'s to be invertible,
then it is possible to get more than one partially ordered set, but clearly
the number is finite and is bounded above by $3^{n}$.

\begin{lemma}
\label{goodpositionyieldsuniquepo}Let $G$ be a finitely generated group with
finitely generated subgroups $H_{1},\ldots,H_{n}$. For $i=1,\ldots,n$, let
$X_{i}$ be a nontrivial $H_{i}$--almost invariant subset of $G$, and suppose
that, for distinct $i$ and $j$, the $G$--orbits of $X_{i}$ and $X_{j}$ are not
parallel. Suppose that $X_{i}$ is equivalent to $Y_{i}$ and to $Z_{i}$ such
that the $Y_{i}$'s are in good position and the $Z_{i}$'s are in good
position. Further suppose that, for each $i$, $Y_{i}$ and $Z_{i}$ are not
invertible. Then there is a $G$--equivariant bijection $\varphi:E(Y_{1}%
,\ldots,Y_{n})\rightarrow E(Z_{1},\ldots,Z_{n})$ which preserves the partial
order $\leq$ and preserves equivalence classes.
\end{lemma}

\begin{proof}
As discussed just after Definition \ref{defnofinvertible}, we can define a
$G$--equivariant bijection $\varphi$ from $E(Y_{1},\ldots,Y_{n})$ to
$E(Z_{1},\ldots,Z_{n})$ by sending $gY_{i}$ to $gZ_{i}$, and $gY_{i}^{\ast}$
to $gZ_{i}^{\ast}$, for each $i$ and for every $g$ in $G$, and $\varphi$ also
preserves complementation and equivalence classes. The proof of Lemma
\ref{goodpositionyieldsuniquepoifn=1} shows how to modify $\varphi$ to be
order preserving when restricted to each $E(Y_{i})$. If $\varphi$ is not
itself order preserving, there are elements $U\ $and $V$ of $E(Y_{1}%
,\ldots,Y_{n})$ such that $U\leq V$ but $\varphi U\nleq\varphi V$. As
$\varphi$ preserves equivalence classes, this implies that the pair $(U,V)$
has two small corners, so that $U$ is equivalent to $V$ or $V^{\ast}$. Let $i$
and $j$ denote those integers such that $U$ is a translate of $Y_{i}$ or
$Y_{i}^{\ast}$ and $V$ is a translate of $Y_{j}$ or $Y_{j}^{\ast}$. If $i=j$,
this contradicts the fact that $\varphi$ is order preserving when restricted
to each $E(Y_{i})$. If $i\neq j$, this contradicts our hypothesis that the
$G$--orbits of $X_{i}$ and $X_{j}$ are not parallel. These contradictions show
that $\varphi$ must be order preserving, as required.
\end{proof}

\section{Constructing cubings from almost invariant sets in good position}

\label{CubingConstruction}As in the previous section, we consider a finitely
generated group $G$ with finitely generated subgroups $H_{1},\ldots,H_{n}$.
For $i=1,\ldots,n$, let $X_{i}$ be a nontrivial $H_{i}$--almost invariant
subset of $G$, and let $E=\{gX_{i},gX_{i}^{\ast}:g\in G,1\leq i\leq n\}$. In
\cite{Sageev-cubings}, Sageev gave a construction of a cubing from $E$, which
we outlined in section \ref{cubings}. A key ingredient of his construction was
the use of the partial order induced by inclusion on $E$. In the previous
section, we established that given a finite family of nontrivial almost
invariant sets, there exists an equivalent family in good position, and, if
the $X_{i}$'s are in good position, we described a new partial order on $E$.
In this section, we describe a variant of Sageev's construction which uses
this new partial order. We will see from the discussion immediately after the
proof of Theorem \ref{thm:Embedding} that this gives a cubing which is minimal
in a natural sense, and in most cases it is canonically associated to the
equivalence classes of the $X_{i}$'s.

Now suppose that the $X_{i}$'s are in good position and consider $E$ with the
partial order of almost inclusion discussed in the previous section. As in
section \ref{cubings}, let $\Lambda^{(0)}$ denote the collection of all
ultrafilters on $E$, defined using the new partial order. Exactly as in
section \ref{cubings}, we can inductively construct a cubed complex $\Lambda$
whose vertex set is $\Lambda^{(0)}$. Again $\Lambda$ will not be connected,
but we wish to pick out a component $L$ which corresponds in a natural way to
the component $C$ picked out in the previous case. In fact the vertices of
$L,$ like the vertices of $C$, will be characterised as ultrafilters on $E$
which satisfy the descending chain condition. We cannot proceed exactly as
before because the set $V_{g}=\{A\in E:g\in A\}$ need not be an ultrafilter
with respect to the new partial order. For example, it is quite possible that
$g\in A\leq B$, but that $g\notin B$. We will thus need to adjust the
construction of basic vertices.

We will need the following technical lemma, which will allow us to start by
constructing an ultrafilter for all but a finite number of elements of $E$.

\begin{lemma}
\label{lem:ContainmentNeighborhood} There exists $R>0$ such that if $A,B\in E$
and $A\leq B$ and if $g\in A$ such that $N_{R}(g)\subset A$, then $g\in B$.
\end{lemma}

\begin{proof}
As $A\leq B$, we also have $B^{\ast}\leq A^{\ast}$. Now Lemma
\ref{containedinboundednbhdofA} tells us that there is $D>0$ such that
$B^{\ast}\subset N_{D}(A^{\ast})=A^{\ast}\cup N_{D}(\delta A)$. If $g$ lies in
$A$ but not in $B$, it follows that $g$ lies in $N_{D}(\delta A)$. This
implies there is a point $h$ of $A^{\ast}$ such that $d(g,h)\leq D+1$, so that
$N_{D+1}(g)$ is not contained in $A$. Thus the lemma holds with $R=D+1$.
\end{proof}

We are now ready to describe the special ultrafilters which will pick out the
component $L$ of $\Lambda$ which corresponds to $C$. Given $g\in G$, we want
to describe an ultrafilter $W_{g}$ which will be almost the same as the set
$V_{g}=\{A\in E:g\in A\}$. Consider first the ball $N=N_{R}(g)$ of radius $R$
about $g$ in the Cayley graph of $G$, where $R$ is as in Lemma
\ref{lem:ContainmentNeighborhood} above. We let%
\[
E_{R}=\{A\in E|{\delta}A\cap N\not =\emptyset\}.
\]

We then denote $E-E_{R}$ by $E_{R}^{\ast}$. As $E$ consists of the translates
of a finite family of $X_{i}$'s and their complements, it follows that $E_{R}$
is finite.

Now for each pair $\{A,A^{\ast}\}$ of elements of $E$ we need to decide
whether or not $A$ or $A^{\ast}$ is in $W_{g}$, consistent with condition 2)
of Definition \ref{defnofultrafilter}. We will make this decision first for
pairs $(A,A^{\ast})$ in $E_{R}^{\ast}$. As in the definition of $V_{g}$, we do
this by taking those elements that contain $g$. That is, let%
\[
U_{g}=\{A\in E_{R}^{\ast}|g\in A\}.
\]

Note that if $A\in U_{g}$, then $N_{R}(g)\subset A$.

\begin{lemma}
\medskip$U_{g}$ is an ultrafilter on $E_{R}^{\ast}$.
\end{lemma}

\begin{proof}
For each pair $\{A,A^{\ast}\}\in E_{R}^{\ast}$, we either have $g\in A$ or
$g\in A^{\ast}$, so that condition 1) of Definition \ref{defnofultrafilter}
holds. Now suppose that $A\in U_{g}$, $B\in E_{R}^{\ast}$ and $A\leq B$. Then
$N_{R}(g)\subset A$, so that Lemma \ref{lem:ContainmentNeighborhood} tells us
that $g\in B$. Hence $B\in U_{g}$, and we have shown that condition 2) of
Definition \ref{defnofultrafilter} holds.
\end{proof}

We now wish to complete $U_{g}$ to an ultrafilter $W_{g}$ on all of $E$. There
are only finitely many pairs $\{A,A^{\ast}\}$ about which we need to make a
decision as to whether $A$ or $A^{\ast}$ is in $W_{g}$.

First of all, for each $B\in E_{R}$ for which there exists $A\in U_{g}$, with
$A\leq B$, we add $B$ to $U_{g}$. That is, set%
\[
U_{1}=U_{g}\cup\{B\in E_{R}\ |\ \exists\ A\in U_{g},A\leq B\}.
\]

\begin{lemma}
$U_{1}$ is an ultrafilter on the set $U_{1}\cup U_{1}^{\ast}$, where
$U_{1}^{\ast}$ denotes the set $\{X^{\ast}:X\in U_{1}\}$.
\end{lemma}

\begin{proof}
By construction $U_{1}$ satisfies condition 2) of Definition
\ref{defnofultrafilter}, namely that if $A\in U_{1}$ and $A\leq B$ then $B\in
U_{1}$. We claim that $U_{1}$ also satisfies condition 1) of Definition
\ref{defnofultrafilter}, namely that we do not have $B\in U_{1}$ and $B^{\ast
}\in U_{1}$. For if this occurs, we have $A_{1}$ and $A_{2}$ in $U_{g}$, with
$A_{1}\leq B$ and $A_{2}\leq B^{\ast}$. Thus we have $A_{1}\leq B\leq
A_{2}^{\ast}$. As $N_{R}(g)\subset A_{1}$, Lemma
\ref{lem:ContainmentNeighborhood} tells us that $g\in A_{2}^{\ast}$, which
contradicts the fact that $g\in A_{2}$. It follows that $U_{1}$ is an
ultrafilter on $U_{1}\cup U_{1}^{\ast}$, as required.
\end{proof}

Now let $V_{1}$ denote the collection of the remaining elements of $E$, so
that $V_{1}=E-(U_{1}\cup U_{1}^{\ast})$, and let $A_{1}$ denote a minimal
element of $V_{1}$. We form $U_{2}$ by adding $A_{1}$ to $U_{1}$ and then
adding every $B\in V_{1}$ such that $A_{1}\leq B$.

\begin{lemma}
$U_{2}$ is an ultrafilter on the set $U_{2}\cup U_{2}^{\ast}$.
\end{lemma}

\begin{proof}
Clearly $U_{2}$ does not contain $B$ and $B^{\ast}$, for any $B$ in $U_{2}\cup
U_{2}^{\ast}$, and so $U_{2}$ satisfies condition 1) of Definition
\ref{defnofultrafilter}. We will show that it also satisfies condition 2). For
suppose $C\in U_{2}$ and $C\leq D$, where $D\in U_{2}\cup U_{2}^{\ast}$. If
$C\in U_{1}$, then the definition of $U_{1}$ implies that $D\in U_{1}$ also
and hence $D\in U_{2}$. If $C\notin U_{1}$, and $D\notin U_{1}\cup U_{1}%
^{\ast}$, then $D\in U_{2}$ by our construction. If $C\notin U_{1}$ and $D\in
U_{1}^{\ast}$, then $D^{\ast}\leq C^{\ast}$ and $D^{\ast}\in U_{1}$, which
implies that $C^{\ast}\in U_{1}$. Thus $C^{\ast}\in U_{2}$ which contradicts
our assumption that $C\in U_{2}$.
\end{proof}

Next let $V_{2}$ denote the collection of the remaining elements of $E$, so
that $V_{2}=E-(U_{2}\cup U_{2}^{\ast})$, and let $A_{2}$ denote a minimal
element of $V_{2}$. We form $U_{3}$ by adding $A_{2}$ to $U_{2}$ and then
adding every $B\in V_{2}$ such that $A_{2}\leq B$. As above, $U_{3}$ is an
ultrafilter on the set $U_{3}\cup U_{3}^{\ast}$.

We continue in this way until all the elements of $E$ have been exhausted. The
resulting subset $W_{g}$ of $E$ is then an ultrafilter on $E$.

Note that $W_{g}$ is not determined by $g$. The construction of $U_{2}$ and
its successors involves making choices of minimal elements. Thus, for each $g$
in $G$, the above construction will yield finitely many such ultrafilters
$W_{g}$. A vertex $W_{g}$ of $\Lambda$ constructed in this way is called a
\textit{basic vertex}. As one sees from the construction, it agrees with the
notion of a basic vertex in the original construction of the cubing in
\cite{Sageev-cubings} except on a finite subset of $E$. The natural action of
$G$ on $E$ preserves the partial order of almost inclusion, and so induces an
action of $G$ on $\Lambda$.

Next we need to show that the basic vertices of $\Lambda$ all lie in a single
component $L$. Recall that any two basic vertices of the cubed complex $K$
constructed by Sageev in \cite{Sageev-cubings} agree except on a finite number
of pairs of elements of $E$. Now each basic vertex $W_{g}$ of $\Lambda$
associated to an element $g$ of $G$ by the above construction agrees with the
basic vertex $V_{g}$ of $K$ except on a finite number of pairs of elements of
$E$. It follows that any two basic vertices of $\Lambda$ are ultrafilters on
$(E,\leq)$ which agree except on a finite number of pairs of elements of $E$.
Suppose that $v$ and $v^{\prime}$ disagree on $k$ pairs of elements of $E$.
Then, as discussed in section \ref{cubings}, there is a path of length $k$ in
$\Lambda$ which joins $v$ to $v^{\prime}$. It follows that the basic vertices
of $\Lambda$ all lie in a single component $L$, as required.

Finally one needs to show that $L$ is simply connected and $CAT(0)$. The
argument here is essentially the same as in \cite{Sageev-cubings} and will be
left to the reader.

Having constructed $L$, we want to compare it with the cubing $C$ constructed
by Sageev in \cite{Sageev-cubings}. The first step is the following result.

\begin{theorem}
\label{thm:Embedding}Let $G$ be a finitely generated group with finitely
generated subgroups $H_{1},\ldots,H_{n}$. For $i=1,\ldots,n$, let $X_{i}$ be a
nontrivial $H_{i}$--almost invariant subset of $G$, and let $E=\{gX_{i}%
,gX_{i}^{\ast}:g\in G,1\leq i\leq n\}$. Suppose that the $X_{i}$'s are in good
position. Let $(E,\subset)$ denote the set $E$ with the partial order given by
inclusion, and let $(E,\leq)$ denote the set $E$ with the partial order given
by almost inclusion, as described in section \ref{CubingConstruction}. Let $C$
denote the cubing constructed from the poset $(E,\subset)$ as in Sageev's
original construction in \cite{Sageev-cubings}, and let $L$ be the cubing
constructed from the poset $(E,\leq)$ as in the previous section. Then there
is a natural $G$--equivariant embedding $L\rightarrow C$.
\end{theorem}

\begin{proof}
Let $K$ denote the cubed complex constructed from $(E,\subset)$, and let
$\Lambda$ denote the cubed complex constructed from $(E,\leq)$, so that $C$ is
a component of $K$ and $L$ is a component of $\Lambda$. We claim first that a
vertex of $\Lambda$ is a vertex of $K$. For if $V$ is an ultrafilter on
$(E,\leq)$, then $V$ is a subset of $E$ which satisfies the following conditions,

\begin{itemize}
\item For any $A\in E$ either $A\in V$ or $A^{\ast}\in V$, but not both.

\item If $A\in V$ and $A\leq B$, then $B\in V$.
\end{itemize}

Now if $A\subset B$, then certainly $A\leq B$, so it follows immediately that
$V$ is also an ultrafilter on $(E,\subset)$. Thus $\Lambda^{(0)}\subset
K^{(0)}$. The description of the construction of the cubed complexes $K$ and
$\Lambda$ from their vertices shows that this inclusion naturally extends to
an embedding of $\Lambda$ in $K$, and that this embedding is $G$--equivariant.
As any basic vertex of $L$ differs from some basic vertex of $C$ by only
finitely many elements, it follows that they can be joined by a path in $K$.
Thus the embedding of $\Lambda$ in $K$ induces an embedding of $L$ in $C$, as required.
\end{proof}

Note that if we are given two collections of good position almost invariant
sets, $Y_{1},\ldots,Y_{n}$ and $Z_{1},\ldots,Z_{n}$ with $Y_{i}$ equivalent to
$Z_{i}$, such that no $Y_{i}$ or $Z_{i}$ is invertible, Lemma
\ref{goodpositionyieldsuniquepo} provides a $G$--equivariant, order preserving
bijection from $E(Y_{1},\ldots,Y_{n})$ to $E(Z_{1},\ldots,Z_{n})$, which
provides a $G$--equivariant isomorphism from $L_{Y}$ to $L_{Z}$. Thus the
cubing constructed from the poset $(E,\leq)$ is determined solely by the
equivalence classes of the almost invariant sets $X_{1},\ldots,X_{n}$.

Now we are ready to compare our new cubing with the old. Suppose we are given
a family of almost invariant subsets $X_{1},\ldots,X_{n}$ of a group $G$, such
that the $X_{i}$'s are in good position. For simplicity we assume further that
no $X_{i}$ is equivalent to an invertible set, and that no two $G$--orbits of
the $X_{i}$'s are parallel. We have just constructed a cubing $L$ which
depends only on the equivalence classes of the $X_{i}$'s. If we consider
almost invariant subsets $Y_{1},\ldots,Y_{n}$ such that $Y_{i}$ is equivalent
to $X_{i}$, we also have Sageev's original cubing $C(Y_{1},\ldots,Y_{n})$,
which we denote by $C(Y)$ for brevity, and Theorem \ref{thm:Embedding} shows
that $L$ embeds in $C(Y)$ for any choices of the $Y_{i}$'s. Thus $L$ is in
some natural sense smaller than any of the $C(Y)$'s. It is clear that $L$ will
equal $C(Y)$ if the partial orders on $E(Y_{1},\ldots,Y_{n})$ induced by
inclusion and by $\leq$ are the same. This is the same condition as the
$Y_{i}$'s being in very good position as defined at the start of the next
section. Now Lemma \ref{verygoodposition} below states that we can always
choose $Z_{i}$ equivalent to $X_{i}$, so that the $Z_{i}$'s are in very good
position. Thus $L$ equals $C(Z)$ and is a minimal cubing among all the cubings
$C(Y)$ obtained by choosing $Y_{i}$ equivalent to $X_{i}$.

As a simple example, consider the special case discussed in the introduction
where $G$ is the fundamental group of a closed orientable surface $M$, and
$S_{1},\ldots,S_{n}$ are a family $\mathcal{F}$ of disjoint simple closed
curves on $M$, such that no two of the $S_{i}$'s are parallel. We can
associate an almost invariant subset $X_{i}$ of $G$ to $S_{i}$ as described
just before Lemma \ref{verygoodposition}, and we can use $X_{1},\ldots,X_{n}$
to construct cubings $L$ and $C(X)$. In this case, $L$ equals $C(X)$ and is
the dual tree to $\widetilde{\mathcal{F}}$ in $\widetilde{M}$. If we now
homotop the curves to meet each other, the associated almost invariant subsets
$Y_{1},\ldots,Y_{n}$ are equivalent to $X_{1},\ldots,X_{n}$ respectively and
yield a new cubing $C(Y)$ which may no longer be $1$--dimensional. In fact, we
can make this cubing have as high a dimension as we please by homotoping the
$S_{i}$'s to meet in a suitably complicated way.

\begin{remark}
If we define distance functions on $L$ and $C$, by assigning length $1$ to
each edge, then the inclusion of the cubing $L$ into the cubing $C$ is
isometric. For if $v$ and $w$ are two vertices in $L$, then the number of
edges in any $C$--geodesic between $v$ and $w$ equals the number of
hyperplanes of $C$ which separate $v$ from $w$. Similarly the number of edges
in any $L$--geodesic between $v$ and $w$ equals the number of hyperplanes of
$L$ which separate $v$ from $w$. These numbers are equal because the vertices
are ultrafilters and, in both cases, the number of hyperplanes separating the
vertices measures the number of sets in $E$ which need to be replaced by their complements.
\end{remark}

\section{Applications\label{applications}}

We saw in section \ref{newpartialorder} that given a family of almost
invariant sets $X_{1}\ldots,X_{n}$, there is a family of almost invariant sets
$Y_{1}\ldots,Y_{n}$, such that $Y_{i}$ is equivalent to $X_{i}$, and the
$Y_{i}$'s are in good position. This means that if two elements of
$E(Y)=E(Y_{1},\ldots,Y_{n})$ have two of their four corners small, then one is
empty. Thus two elements of $E(Y)$ must cross, be nested or have only one
small corner. In this section, we will show that the third possibility can be
removed. Precisely, we say that the $X_{i}$'s are in \textit{very good
position} if given two elements of $E(X)$, either they cross or they are
nested. This means that the partial orders on $E(X_{1},\ldots,X_{n})$ induced
by inclusion and by $\leq$ are the same. We will show that we can always
arrange this situation by replacing each $X_{i}$ by an equivalent almost
invariant set $Z_{i}$.

As we stated in the introduction, very good position for almost invariant sets
is closely analogous to the properties enjoyed by shortest curves on surfaces
or by least area surfaces in $3$--manifolds. For simplicity, we will discuss
only curves on surfaces. In order to explain the analogy, we first need to
recall how curves on a surface are related to almost invariant sets. Let $F$
denote a closed orientable surface and let $S$ denote a simple closed curve on
$F$. Let $H$ denote the infinite cyclic subgroup of $G=\pi_{1}(F)$ carried by
$S$, and let $F_{H}$ denote the cover of $F$ whose fundamental group is $H$.
Thus $S$ lifts to a circle in $F_{H}$ which we also denote by $S$. Pick a
generating set for $G$ and represent it by a bouquet of circles embedded in
$F$. We will assume that the base point of the bouquet does not lie on $S$.
The pre-image of this bouquet in the universal cover $\widetilde{F}$ of $F$
will be a copy of the Cayley graph $\Gamma$ of $G$ with respect to the chosen
generating set. The pre-image in $F_{H}$ of the bouquet will be a copy of the
graph $H\backslash\Gamma$, the quotient of $\Gamma$ by the action of $H$ on
the left. Consider the closed curve $S$ on $F_{H}$. Let $P$ denote the set of
all vertices of $H\backslash\Gamma$ which lie on one side of $S$. Then $P$ has
finite coboundary, as $\delta P$ equals exactly the edges of $H\backslash
\Gamma$ which cross $S$. Hence $P$ is an almost invariant subset of
$H\backslash G$. Let $X$ denote the pre-image of $P$ in $\Gamma$, so that $X$
equals the set of vertices of $\Gamma$ which lie on one side of the line $l$.
Then $X$ is a $H$--almost invariant subset of $G$. If $S$ is not simple, but
we choose it to be shortest in its homotopy class, its lift to $F_{H}$ will
still be simple, so that the same construction can be made. Now the fact that
$S$ is shortest implies that, for each $g\in G$, the translate $gl$ of the
line $l$ in $\widetilde{F}$ must equal $l$, be disjoint from $l$ or meet $l$
transversely in a single point. If $gl=l$, it follows that the translate $gX$
of $X$ must equal $X$ (it cannot equal $X^{\ast}$ as $F$ is orientable). If
$gl$ is disjoint from $l$, it follows that $gX$ and $X$ are nested. If $gl$
meets $l$ transversely in a single point, it follows that $X$ and $gX$ cross
each other. We conclude that if $S$ is shortest in its homotopy class, then
$X$ is in very good position.

\begin{lemma}
\label{verygoodposition}Let $G$ be a finitely generated group with finitely
generated subgroups $H_{1},\ldots,H_{n}$. For $i=1,\ldots,n$, and let $X_{i}$
be a nontrivial $H_{i}$--almost invariant subset of $G$. Then, for each $i$,
there exists a $K_{i}$--almost invariant subset $Z_{i}$ of $G$ which is
equivalent to $X_{i}$, such that the $Z_{i}$'s are in very good position.
\end{lemma}

\begin{proof}
For simplicity we will consider the case when $n=1$, and will denote $X_{1}$
by $X$ and $H_{1}$ by $H$. The general case is essentially the same. Start
with $Y$ in good position such that $Y$ is equivalent to $X$. Then construct
the cubing $C$ given by $Y$ and using the poset $(E(Y),\subset)$. As discussed
just before Lemma \ref{lem:ConstructinAISets}, there is a hyperplane
$\mathcal{H}$ of $C$ and a half-space $\mathcal{H}^{+}$ determined by
$\mathcal{H}$ such that a vertex of $C$ lies in $\mathcal{H}^{+}$ if and only
if, when regarded as an ultrafilter on $E$, it contains $Y$. Further, for any
vertex $v$ of $C$, the set $Y_{v}=\{g\in G\mid g(v)\in\mathcal{H}^{+}\}$ is
$H$--almost invariant and equivalent to $Y$.

Next consider the cubing $L$ given by $Y$ and using the poset $(E(Y),\leq)$,
as constructed in section \ref{CubingConstruction}. Recall that $\mathcal{H}$
is associated to an equivalence class $F$ of edges of $C$ given by the
equivalence relation generated by saying that two edges are equivalent if they
are opposite edges of a square in $C$. Now two edges of $L$ are opposite edges
of a square in $L$ if and only if they are opposite edges of a square in $C$.
It follows that if $f$ is an edge of $F$ which also lies in $L$, then the
equivalence class of $f$ in $L$ is precisely $F\cap L$. Let $\mathcal{K}$
denote the hyperplane in $L$ associated to this equivalence class. Then it
follows that $\mathcal{H}^{+}\cap L$ equals one of the two half-spaces in $L$
determined by $\mathcal{K}$. We denote this half-space by $\mathcal{K}^{+}$.
Pick a vertex $w$ of $L$, and apply Lemma \ref{lem:ConstructinAISets} to
obtain a new almost invariant set $Z$ over $H$ equal to $Y_{w}=\{g\in G\mid
g(w)\in\mathcal{K}^{+}\}$. Since the inclusion of $L$ in $C$ is $G$%
--equivariant, $Z$ may also be viewed as the set $\{g\in G\mid g(w)\in
\mathcal{H}^{+}\}$. Now Lemma \ref{lem:ConstructinAISets} tells us that
$Y_{v}$ and $Y_{w}$ are equivalent. As $Y_{v}$ is equivalent to $Y$, and
$Z=Y_{w}$, it follows that $Z$ is equivalent to $Y$.

In particular, two translates of $Y$ in $G$ are almost nested if and only if
the corresponding translates of $Z$ are almost nested. However, we claim that
if two translates of $Y$ in $G$ are almost nested then the corresponding
translates of $Z$ are actually nested. This means exactly that $Z$ is in very
good position. To prove our claim, suppose, for example, that $aY\leq Y$. We
need to show that $aZ\subset Z$. Recall that $Z$ can be viewed as $\{g\in
G\mid g(w)\in\mathcal{H}^{+}\}$. The description of the vertices of
$\mathcal{H}^{+}$ given in the first paragraph of this proof shows that
$Z=\{g\in G\mid Y\in g(w)\}$. Thus
\begin{align*}
aZ  &  =\{ag\in G\mid g(w)\in\mathcal{H}^{+}\}=\{g\in G\mid a^{-1}%
g(w)\in\mathcal{H}^{+}\}\\
&  =\{g\in G\mid Y\in a^{-1}g(w)\}=\{g\in G\mid aY\in g(w)\}.
\end{align*}
As $w$ is an ultrafilter on $(E(Y),\leq)$, so is $g(w)$. As $aY\leq Y$, it
follows that if $aY\in g(w)$, then $Y\in g(w)$. Thus $aZ\subset Z\ $as
claimed. As this holds for all $a$, and analogous arguments apply if $aY\leq
Y^{\ast}$, $aY^{\ast}\leq Y$ or $aY^{\ast}\leq Y^{\ast}$, it follows that $Z$
is in very good position, completing the proof of the lemma.
\end{proof}

We next consider applications which strengthen results of Niblo in
\cite{Niblo} on the existence of splittings of a given group. Let $H$ be a
finitely generated subgroup of a finitely generated group $G$, and let $X$ be
a nontrivial $H$--almost invariant subset of $G$. In \cite{Niblo}, Niblo
defined a group $T(X)$ which is the subgroup of $G$ generated by $H$ and
$\{g\in G:gX$ and $X$ are not nested\}. He proved, using Sageev's construction
of cubings, that if $T(X)\neq G$, then $G$ splits over a subgroup of $T(X)$.
One can also define $S(X)$ to be the subgroup of $G$ generated by $H$ and
$\{g\in G:gX$ crosses $X\}$. Clearly $S(X)$ is contained in $T(X)$. Further
they are equal if $X$ is in very good position. Thus the fact that $X$ is
equivalent to an almost invariant set in very good position yields a
strengthening of Niblo's result in which one can replace $T(X)$ by $S(X)$,
i.e. one can replace the condition of not being nested by the condition of
crossing. This strengthening was obtained previously by Scott and Swarup in
\cite{SS-regnbhds} using their theory of regular neighbourhoods, but the
present argument is more elementary.

In \cite{Niblo}, Niblo proved an analogous result for two almost invariant
subsets of a finitely generated group $G$. Again he used Sageev's construction
of cubings. Let $K$ is another finitely generated subgroup of $G$, and let $Y$
be a nontrivial $K$--almost invariant subset of $G$. Suppose that any
translate of $X$ and any translate of $Y$ are nested. Then $G$\ splits over a
subgroup of $T(X)$\ and over a subgroup of $T(Y)$. More precisely $G$ is the
fundamental group of a graph of groups with two edges such that the edge
groups are conjugate into $T(X)\ $and $T(Y)$ respectively. As above, the fact
that $X$ and $Y$ can be replaced by equivalent almost invariant sets in very
good position means that the assumption that any translate of $X$ and any
translate of $Y$ are nested can be replaced by the assumption that any
translate of $X$ and any translate of $Y$ do not cross. This strengthening was
also obtained previously by Scott and Swarup in \cite{SS-regnbhds} using their
theory of regular neighbourhoods.

Finally, we state a result which generalises a result of Dunwoody and Roller
in \cite{D-Roller} and strengthens a result of Niblo in {\cite{Niblo}}.

\begin{theorem}
Let $G$ be a finitely generated group with a finitely generated subgroup $H$
and a nontrivial $H$--almost invariant subset $X$. If $\{g\in G:gX$ crosses
$X\}$ lies in $Comm_{G}(H)$, the commensuriser of $H$ in $G$, then $G$ splits
over a subgroup commensurable with $H$.
\end{theorem}

In \cite{Niblo}, Niblo proved this result on the stronger assumption that
$\{g\in G:gX$ and $X$ are not nested\} lies in $Comm_{G}(H)$. In
\cite{D-Roller}, Dunwoody and Roller proved the special case of this result
when $G$ commensurises $H$. One way to prove the result stated above is simply
to apply Niblo's result using the fact that $X$ is equivalent to an almost
invariant subset of $G$ in very good position. Alternatively, as Niblo's
argument used Sageev's construction of cubings, one could obtain the
strengthened result more directly by using our new cubing in place of Sageev's
in Niblo's argument.

\end{document}